\input amstex
\input amsppt.sty
\magnification\magstep1
\pagewidth{12.9 true cm}
\pageheight{19.0 true cm}
\pageno=1
\NoBlackBoxes
\tolerance=10000
\TagsOnRight

\leftheadtext{S. S. Dragomir, Y. J. Cho and S. S. Kim}
\rightheadtext{The Cauchy-Bunyakovsky-Schwarz inequality}
\topmatter
\title
Some Reverses of the Cauchy-Bunyakovsky-Schwarz inequality in 2-inner
product spaces
\endtitle
\endtopmatter
\author 
Sever S. Dragomir, Yeol Je Cho and Seong Sik Kim
\endauthor

\abstract
In this paper, some reverses of the  Cauchy-Bunyakovsky-Schwarz inequality in 2-inner product spaces are given. Using this framework, some applications for determinantal integral inequalities are also provided.
\endabstract
\footnote""{{2000 AMS Subject Classification:} Primary 46C05, 46C97. Secondary 26D15, 26D10.}
\footnote""{{Key Words and Phrases:} 2-inner product spaces, Cauchy-Bunyakovsky-Schwarz inequality, Determinantal integral inequalities.}
\footnote""{The corresponding author: yjcho\@nongae.gsnu.ac.kr (Y. J. Cho) or sskim\@dongeui.ac.kr (S. S. Kim).}
\endtopmatter
\document
\centerline{\bf 1. Introduction}
\vskip 2mm

The concepts of 2-inner products and 2-inner product spaces have been intensively studied by many authors in the last three decades.
A systematic presentation of the recent results related to the theory of 2-inner product spaces as well as an extensive list of
the related references can be found in the book [3].
\vskip 0.2cm 

We recall here the basic definitions and the elementary properties of 2-inner product spaces that will be used in the sequel (see [4]).
\vskip 2mm

Let $X$ be a linear space of dimension greater than $1$ over the number filed $K$, where $K=R$ or $K=C$. Suppose that
$(\cdot,\cdot|\cdot)$ is a $K$-valued function on $X \times X \times X$ satisfying the following conditions:
\roster
\item"{$(2I_1)$}" ($x,x |z) \ge 0$,
$(x,x |z)=0$ if and only if $x$ and $z$ are linearly dependent,
\item"{$(2I_2)$}" $(x,x|z) = (z,z |x)$,
\item"{$(2I_3)$}" $(x,y |z) = \overline{(y,x | z)}$,
\item"{$(2I_4)$}" $(\alpha x,y |z)= \alpha(x,y|z)$ for any scalar $\alpha \in K$,
\item"{$(2I_5)$}" $(x+x^\prime,y |z) = (x,y|z)+(x^\prime,y |z)$,
\endroster
 where $x, x',y, z \in X$. 

The functional $(\cdot,\cdot|\cdot)$ is called a {\it 2-inner product} and ($X,(\cdot,\cdot|\cdot))$ a {\it 2-inner product space} (or {\it 2-pre-Hilbert space}).
\vskip 0.2cm 

Some basic properties of the 2-inner product spaces are as follows:
\vskip 2mm

(1) If $K=R$, then (2$I_3$) reduces to $(x,y|z) = (y,x|z)$.

(2) From (2$I_3$) and (2$I_4$), we have $(0,y|z)=0, (x,0|y)=0$ and also
$$
(x,\alpha y|z)=\overline \alpha(x,y|z).
\leqno(1.1)
$$

(3) Using (2$I_3$)$\sim$(2$I_5$), we have
$$
(z,z|x\pm y)=(x\pm y,x\pm y|z)=(x,z|z)+(y,y|z) \pm 2Re(x,y|z)
$$
and
$$
Re(x,y|z)={1 \over 4}[(z,z|x+y)-(z,z|x-y)].
\leqno(1.2)
$$

In the real case $K=R$, (1.2) reduces to
$$
(x,y|z)={1 \over 4}[(z,z|x+y)-(z,z|x-y)]
\leqno(1.3)
$$
and, using this formula, it is easy to see that, for any $\alpha \in R$,
$$
(x,y|\alpha z)=\alpha^2(x,y|z).
\leqno(1.4)
$$

In the complex case $K=C$, using (1.1) and (1.2), we have
$$
Im(x,y|z)=Re[-i(x,y|z)]={1 \over 4}[(z,z|x+iy)-(z,z|x-iy)],
$$
which, in combination with (1.2), yields
$$
\aligned
(x,y|z)&={1 \over 4}[(z,z|x+y)-(z,z|x-y)]\\
&\quad +{i \over 4}[(z,z|x+iy)-(z,z|x-iy)].
\endaligned
\leqno(1.5)
$$

Using (1.5) and (1.1), we have, for any $\alpha \in C$,  that
$$
(x,y|\alpha z)=|\alpha|^2(x,y|z).
\leqno(1.6)
$$
However, for any $\alpha \in R$, (1.6) reduces to (1.4). Also, it follows  from (1.6) that
$$
(x,y|0)=0.
$$

(4) For any  given vectors $x,y,z \in X$, consider the vector $u=(y,y|z)x-(x,x|z)y$. By (2$I_1$), we know that $(u,u|z)\ge 0$. It
is obvious that the inequality $(u,u|z) \ge 0$ can be rewritten as
$$
(y,y|z)[(x,x|z)(y,y|z)-|(x,y|z)|^2] \ge 0.
\leqno(1.7)
$$
If $x=z$, then (1.7) becomes
$$
-(y,y|z)|(z,y|z)|^2 \ge 0,
$$
which implies that
$$
(z,y|z)=(y,z|z)=0,
\leqno(1.8)
$$
provided $y$ and $z$ are linearly independent. Obviously, when $y$ and $z$ are linearly dependent, (1.8) holds too.

 Now, if $y$ and $z$ are linearly  independent, then $(y,y|z)>0$ and, from (1.2), it follows the Cauchy-Bunyakovsky-Schwarz inequality (shortly,
the CBS-inequality) for 2-inner products:
$$
|(x,y|z)|^2 \le (x,x|z)(y,y|z).
\leqno(1.9)
$$
Using (1.8), it is easy to see that (1.9) is trivially fulfilled when $y$ and $z$ are linearly dependent. Therefore, the inequality
(1.9) holds for any three vectors $x,y,z \in X$ and is strict unless $u=(y,y|z)x-(x,y|z)y$ and $z$ are linearly dependent. In
fact, we have the equality in (1.9) if and only if the three vectors $x,y$ and $z$ are linearly dependent. 
\vskip 0.2cm 

In any given 2-inner product space $(X,(\cdot,\cdot|\cdot))$, we can define a function $\|\cdot|\cdot\|$ on $X \times X$ by
$$
\|x|z\|=\sqrt{(x,x|z)}
\leqno(1.10)
$$
for all $x,z \in X$.  It is easy to see that this function satisfies the following conditions: 
\roster
\item"{$(2N_1)$}" $\| x|z\| = 0$ if and only if $x$ and $z$ are linearly dependent,
\item"{$(2N_2)$}" $\|x|z\|=\|z|x\|$, 
\item"{$(2N_3)$}" $\|\alpha x|z\|=|\alpha| \|x|z\|$ for all scalar $\alpha \in K$,
\item"{$(2N_4)$}" $\|x+x'|z\| \leq \|x|z\|+\|x'|z\|$. 
\endroster

Any function $\|\cdot,\cdot\|$ satisfying the conditions (2N$_1$)$\sim$(2N$_4$) is called a {\it 2-norm} on $X$ and $(X,\|\cdot|\cdot\|)$ a {\it  linear 2-normed space}. 
\vskip 0.2cm 

For a systematic presentation of the recent results related to the theory of linear 2-normed spaces, see the book [8].
\vskip 2mm

In terms of the 2-norms, the (CBS)-inequality (1.9) can be written as
$$
|(x,y|z)|^2 \le \|x|z\|^2\|y|z\|^2.
\leqno(1.11)
$$

The equality holds in (1.11) if and only if $x,y$ and $z$ are linearly dependent.
\vskip 2mm

 For recent inequalities in 2-inner product spaces, see the papers [1, 2], [4--7], [9--12] and the references therein.
\bigskip

\centerline{\bf 2. Revereses of the CBS-Inequality} 
\vskip 0.2cm

The following reverse of the CBS-inequality in 2-inner product spaces holds.
\vskip 2mm

{\bf Theorem 2.1.} {\it Let $A, a \in K(K=C, R)$ and $x, y,z \in X$, where, as above, $(X, (\cdot,\cdot|\cdot))$ is a 2-inner
product space over $K$. If
$$ 
Re(Ay-x, x-ay|z) \ge 0 
\leqno(2.1)
$$
or, equivalently,
$$
\Big\|x-{a+A \over 2}y\Big|z\Big\| \le {1 \over 2}|A-a|\|y|z\|
\leqno(2.2)
$$
holds, then we have the inequality
$$
0 \le \|x|z\|^2\|y|z\|^2-|(x,y|z)|^2 \le {1 \over 4}|A-a|^2\|y|z\|^4.
\leqno(2.3)
$$
The constant ${1\over 4}$ is sharp in $(2.3)$ in the sense that it cannot be replaced by a smaller constant.}
\vskip 0.2cm 

{\bf Proof.} Consider the vectors $x, u, U, z \in X$. We observe that
$$
\split {1\over 4}
&\|U-u|z\|^2-\Big\|x-{U+u \over 2}\Big|z\Big\|^2\\
&=Re(U-x,x-u|z)\\
&= Re[(x, u|z)+(U, x|z)] - Re(U, u|z) -\|x| z\|^2.
\endsplit
$$
Therefore, $Re(U-x, z-u|z) \ge 0$ if and only if 
$$
\Big\|x-{U+u \over2}\Big|z\Big\| \le {1\over 2}\|U-u|z\|.
$$ 
If we apply this to the vectors $U= Ay$ and $u =ay$, then  we deduce that the inequalities (2.1) and (2.2) are equivalent, as stated.

 Now, let us consider the real numbers
 $$ I_1=Re[(A\|y|z\|^2-(x,y|z))(\overline{(x,y|z)}-\overline a\|y|z\|^2)]$$
 and
     $$ 
 I_2 = \|y|z\|^2 Re(Ay-x, x-ay|z).
$$

A simple calculation shows that
  $$ 
I_1=\|y|z\|^2 Re[A(x,y|z)+\overline a(x,y|z)] -|(x,y|z)|^2-\|y|z\|^4Re(A \overline a)
$$
and
$$
\aligned
 I_2&=\|y|z\|^2 Re[A \overline{(x,y|z)}+\overline a(x,y|z)]\\
&\quad -\|x|z\|^2\|y,z\|^2-\|y,z\|^4 Re(A \overline a),
\endaligned
$$
which give
$$
I_1 -I_2=\|x|z\|^2\|y|z\|^2-|(x,y|z)|^2
\leqno(2.4)
$$
for any $x, y, z \in X$ and $a, A \in K$. If (2.2) holds, then $I_2 \ge 0$ and thus
$$
\aligned
&\|x|z\|^2\|y|z\|^2-|(x, y|z)|^2 \\
&\le Re[(A\|y|z\|^2-(x,y|z))(\overline{(x,y|z)}-\overline a\|y|z\|^2).
\endaligned
\leqno(2.5)
$$

If we use the elementary inequality $Re(\alpha \overline \beta) \le {1 \over 4}|\alpha+\beta|^2$ for any $\alpha,\beta \in K(K=R,C)$,
then we have that
$$
\aligned
&Re[(A\|y|z\|^2-(x,y|z))(\overline{(x,y|z)}-\overline a\|y,z\|^2)]\\
&\le {1 \over 4}|A-a|^2\|y|z\|^4.
\endaligned
\leqno(2.6)
$$

Making use of the inequalities (2.5) and (2.6), we deduce the desired inequality (2.3).

To prove the sharpness of the constant ${1 \over 4}$, assume that
(2.4) holds with a constant $C>0$, i.e.,
$$
\|x|z\|^2\|y|z\|^2 -|(x, y|z)|^2 \le C|A-a|^2\|y|z\|^4,
\leqno(2.7)
$$
where $x, y, z, A$ and $a$ satisfy the hypothesis of the theorem.

Consider $y\in X$ with $\|y|z\|=1$, $a \ne A$, $m \in X$ with $\|m|z\|=1$ and $(y,m|z)=0$ and define
$$
x = {A+a \over 2}y+{A-a \over 2}m.
$$
Then we have
$$
(Ay-x, x-ay|z) = {|A-a|^2 \over 4}(y-m, y+m|z) = 0
$$
 and then the condition (2.1) is fulfilled. From (2.7), we deduce
$$
\aligned
&\Big\|{A+a \over 2}y+{A-a \over 2}m\biggl|z\Big\|^2-\Big|\Big({A+a \over 2}y+{A-a \over 2}m,y\biggl|z\Big)\biggr|^2\\
& \le C|A-a|^2
\endaligned
\leqno(2.8)
$$
and, since
$$
\Big\|{A+a \over 2}y+{A-a \over 2}m\Big|z\Big\|^2=\biggl|{A+a \over 2}\biggr|^2+\biggl|{A-a \over 2}\biggr|^2
$$
and
$$
\biggl|\Big({A+a \over 2}y+{A-a \over 2}m,y\Big|z\Big)\biggr|^2=\biggl|{A+a \over 2}\biggr|^2,
$$
by (2.8),  we have
$$
\biggl|{A-a \over 2}\biggr|^2 \le C|A-a|^2
$$
for any $A,a \in K$ with $a\ne A$, which implies $C \ge {1 \over 4}$. This completes the proof.
\vskip 0.2cm 

Another reverse of the  CBS-inequality in 2-inner product spaces is incorporated in the following theorem:
\vskip 2mm

{\bf Theorem 2.2.} {\it Assume that $x,y,z, a$ and $A$ are the same as in Theorem 2.1. If $Re(\overline a,A)>0$,  then we have the inequalities
$$
\aligned
\|x|z\|\|y|z\|& \le {1 \over 2}\cdot {Re[(\overline A +\overline a)(x,y|z)] \over [Re(\overline a A)]^{\frac{1}{2}}} \\
&\le {1 \over 2}\cdot {|A+a|\over [Re(\overline a A)]^{\frac{1}{2}}}|(x,y|z)|. 
\endaligned
\leqno(2.9)
$$
The constant ${1 \over 2}$ is best possible in both inequalities in the sense that it cannot be replaced by a smaller constant.}
\vskip 2mm

{\bf Proof.} Define
$$
\split
I&=Re(Ay-x,x-a y|z)\\
&=Re[A\overline{(x,y|z)}+\overline a (x,y|z)]-\|x|z\|^2-Re(\overline a A)\|y|z\|^2.
\endsplit
$$
We know that, for a complex number $\alpha \in C, Re(\alpha)=Re(\overline \alpha)$ and thus
$$Re[A\ \overline{(x,y|z)}] =Re[\overline A (x,y|z)],$$
which implies
$$
I=Re[(\overline A+\overline a)(x,y|z)]-\|x|z\|^2 -Re(\overline a A)\|y|z\|^2.
\leqno(2.10)
$$
Since $x,y,z,\alpha,A$ are assumed to satisfy the condition (2.1), by (2.10), we deduce the inequality
$$
\|x|z\|^2+Re(\overline a A)\|y|z\|^2 \le Re[(\overline A+\overline a)(x,y|z)],
$$
which gives
$$ 
\aligned
&{1 \over  [Re(\overline a A)]^{\frac{1}{2}}}\|x|z\|^2+ [Re(\overline a A)]^{\frac{1}{2}}\|y|z\|^2\\
&\le {Re[(\overline A + \overline a)(x,y|z)] \over [Re(\overline aA)]^{\frac{1}{2}}}
\endaligned
\leqno(2.11)
$$
since $Re(\overline a A)>0$. 

On the other hand, by the elementary inequality
$$
\alpha p^2 +{1 \over \alpha}q^2 \ge 2pq
$$
for $p, q \ge 0$ and $\alpha>0$, we have
$$
2\|x|z\|\|y|z\| \le {1 \over [Re(\overline aA)]^{1 \over 2}}\|x|z\|^2
   +[Re(\overline aA)]^{1 \over 2}\|y|z\|^2. 
\leqno(2.12)
$$
Using (2.11) and (2.12), we deduce the first inequality in (2.9). The last part is obvious by the fact that, for $z \in C, |Re(z)| \le |z|$.

To prove the sharpness of the constant ${1 \over 2}$ in the first inequality in (2.9), we assume that (2.9) holds with a constant
$C>0$, i.e.,
$$
\|x|z\|\|y|z\| \le C {{Re[(\overline A + \overline a)(x,y|z)]} \over [Re(\overline aA)]^{1 \over 2}},
\leqno(2.13)
$$
provided $x,y,z,a$ and $A$ satisfy (2.1). If we choose $a=A=1$, $y=x \ne 0$, then obviously (2.1) holds and, from (2.13), we obtain
$$
\|x|z\|^2 \le 2C\|x|z\|^2
\leqno(2.14)
$$
for any linearly independent vectors $x,z \in X$, which implies $C \ge {1 \over 2}$. This completes the proof.
\vskip 2mm

When the constants involved are assumed to be positive, then we may state the following result:
\vskip 2mm

{\bf Corollary 2.3.} {\it Let $M \ge m >0$ and assume that, for $x,y,z \in X$, we have
$$Re(My-x,x-my|z) \ge 0$$
or, equivalently,
$$
\Big\|x-{m+M \over 2}\Big|z\Big\| \le {1 \over 2}(M-m)\|y,z\|.
$$
Then we have the following reverse of the  CBS-inequality
$$
\|x|z\| \|y|z\| \le {1 \over 2}\cdot {M+m \over \sqrt{mM}} Re(x,y|z)
\le {1 \over 2}\cdot {M+m \over \sqrt{mM}} |(x,y|z)|. 
\leqno(2.15)
$$ 
The constant ${1 \over 2}$ is sharp in $(2.15)$.}
\vskip 0.2cm 

Some additive versions of the above are obtained in the following:
\vskip 2mm

{\bf Corollary 2.4.} {\it With the assumptions of  Theorem 2.2, we have the inequalities
$$
\aligned
0 &\le \|x|z\|^2 \|y|z\|^2-|(x,y|z)|^2 \\
&\le {1 \over 4}\cdot {|A-a|^2 \over Re(\overline a A)}|(x,y|z)|^2.
\endaligned
\leqno(2.16)
$$
The constant ${1 \over 4}$  is best possible in $(2.16)$.}
\vskip 2mm

{\bf Corollary 2.5.} {\it With the assumptions of Corollary 2.3, we have
$$
\aligned
0 &\le \|x|z\| \|y|z\|-|(x,y|z)| \le  \|x|z\|\|y|z\|-Re(x,y|z)\\
& \le {1 \over 2}\cdot{(\sqrt M-\sqrt m)^2 \over \sqrt{mM}} Re(x,y|z)
   \le {1 \over 2}\cdot {(\sqrt M-\sqrt m)^2 \over \sqrt{mM}} |(x,y|z)|
\endaligned
\leqno(2.17)
$$
and
$$
\aligned
0 &\le \|x|z\|^2 \|y|z\|^2-|(x,y|z)|^2 \le  \|x|z\|^2\|y|z\|^2-[Re(x,y|z)]^2\\
& \le {1 \over 4}\cdot {(M-m)^2 \over mM} [Re(x,y|z)]^2 \le {1 \over 4}\cdot{(M-m)^2 \over {mM}} |(x,y|z)|^2.
\endaligned
\leqno(2.18)
$$
The constant  ${1 \over 2}$ in $(2.17)$ and the constant  ${1 \over 4}$ in $(2.18)$ are best possible.}
\vskip 2mm

 The third inequality in (2.17) may be used to point out a reverse of the triangle inequality in 2-inner product spaces.
\vskip 2mm

{\bf Corollary 2.6.} {\it Assume that $x,y,z,m$ and $M$ are the same as in Corollary 2.3.  Then we have the following reverse of the triangle
inequality
$$ 
0 \le \|x|z\| + \|y|z\| -\|x+y|z\| \le {{\sqrt M-\sqrt m} \over (mM)^{\frac{1}{4}}}\sqrt{ Re(x,y|z)}.
\leqno(2.19)
$$}
\vskip 2mm

{\bf Proof.} It is easy to see that 
$$0 \le (\|x|z\|+ \|y|z\|)^2-\|x+y|z\|^2 =2[\|x|z\|\|y\|-Re(x,y|z)]
\leqno(2.20)
$$ 
for any $x,y,z \in X$. If the assumptions of Corollary 2.3 hold, then (2.17) is valid and, by (2.20), we deduce
$$
0 \le (\|x|z\|+ \|y|z\|)^2 -\|x+y|z\|^2 \le {({\sqrt M-\sqrt m})^2 \over \sqrt{mM}}{ Re(x,y|z)},
$$
which gives
$$
(\|x|z\|+\|y|z\|)^2 \le \|x+y|z\|^2 + {({\sqrt M-\sqrt m})^2 \over \sqrt{mM}}{ Re(x,y|z)}.
\leqno(2.21)
$$
Taking the square root in (2.21), we have
$$
\split \|x|z\|+\|y|z\|
&\le \sqrt{\|x+y|z\|^2 + {({\sqrt M-\sqrt m})^2 \over \sqrt{mM}}{Re(x,y|z)} }\\
&\le \|x+y|z\| +{{\sqrt M-\sqrt m} \over (mM)^{\frac{1}{4}}}\sqrt{Re(x,y|z)},
\endsplit
$$
from where we deduce the desired inequality (2.21). This completes the proof.
\bigskip 

\centerline{\bf 3. Integral Inequalities}
\vskip 0.2cm 

Let $(\Omega, \Sigma,\mu)$ be a measure space consisting of a set $\Omega$, a $\sigma$-algebra $\Sigma$ of subsets of $\Omega$ and a countably additive and positive measure $\mu$ on $\Sigma$ with valued in $R \cup \{ \infty \}$.  Denote by $L_\phi^2(\Omega)$ the Hilbert space of all real-valued functions
$f$ defined on $\Omega$ that are 2-$\phi$-integrable on $\Omega$,  i.e., $\int_\Omega \phi(s)|\phi(s)|^2 d\mu(s) < \infty$, where $\phi:\Omega \to [0,\infty)$ is a measurable function on $\Omega$. 

We can introduce the following 2-inner product on $L^2_\phi(\Omega)$ 
$$
\aligned
&(f,g|h)_\phi\\
&={1 \over 2}\int_\Omega\int_\Omega\phi(x)\phi(y) \left|\matrix f(x) & f(y)\\
                                                              h(x) & h(y)\\
                                                   \endmatrix\right|
\cdot
        \left|\matrix g(x) & g(y)\\
                       h(x) & h(y)\\
\endmatrix\right|
d\mu(x)d\mu(y), 
\endaligned
\leqno(3.1)
$$
where by $\left|\matrix f(x) & f(y)\\
                   h(x) & h(y)\\
                \endmatrix\right|$
we understand the determinant of the matrix
          $\left[\matrix f(x) & f(y)\\
                   h(x) & h(y)\\
            \endmatrix\right]$,
generating the 2-norm
$$\|f|h\|_\phi =\biggl({1 \over 2}\int_\Omega\int_\Omega\phi(x)\phi(y)
         \left|\matrix f(x) & f(y)\\
                       h(x) & h(y)\\
                \endmatrix \right|^2
d\mu(x)d\mu(y)\biggr)^{1 \over 2}.
\leqno(3.2)
$$

A simple computation with integrals shows that
$$
(f,g|h)_\phi=\left|\matrix \int_\Omega \phi(x)f(x)g(x)d\mu(x)&
                           \int_\Omega \phi(x)f(x)h(x)d\mu(x)\\
                           \int_\Omega \phi(x)g(x)h(x)d\mu(x) &
                           \int_\Omega \phi(x)h^2(x)d\mu(x)\\
             \endmatrix\right|
$$
and
$$
\|f|h\|_\phi=\left|\matrix  \int_\Omega \phi(x)f^2(x)d\mu(x)&
                        \int_\Omega \phi(x)f(x)h(x)d\mu(x)\\
                        \int_\Omega \phi(x)f(x)h(x)d\mu(x) &
                        \int_\Omega \phi(x)h^2(x)d\mu(x)\\
               \endmatrix \right|^{1 \over 2}.
$$

We recall that the pair of function $(q,p) \in L^2_\phi(\Omega)\times L^2_\phi(\Omega)$ is called {\it synchronous} if
$$
(q(x)-q(y))(p(x)-p(y) )\ge 0
$$
for almost every  $x,y \in \Omega$.  If $\Omega =[a,b]$ is an interval of real numbers, then a sufficient condition of synchronicity for $(p,q)$ is that they are monotonic in the same sense, i.e., both of them are increasing or decreasing on $[a,b]$.

Now, suppose that $h \in L^2_\phi(\Omega)$ is such that $h(x) \ne 0$ for almost every $x \in \Omega$. Then, by (3.1), we have
$$
\aligned
&(f,g|h)_\phi\\
&={1 \over 2}\int_\Omega\int_\Omega\phi(x)\phi(y)
         h^2(x) h^2(y)\biggl({f(x) \over h(x)}-{f(y) \over h(y)}\biggr)\\
&\quad \times  \biggl({g(x) \over h(x)}-{g(y) \over h(y)}\biggr)d\mu(x)d\mu(y)
\endaligned
\leqno(3.3)
$$
and thus a sufficient condition for the inequality
$$
(f,g|h)_\phi \ge 0 
\leqno(3.4)
$$
to hold is that the pair of function $({f \over h}, {g \over h})$ to be synchronous.
\vskip 0.2cm 

We are able now to state some integral inequalities that can be derived using the general framework presented above.
\vskip 2mm

{\bf Proposition 3.1.} {\it Let $M>m>0$ and $f,g,h \in L^2_\phi(\Omega)$ so that the functions
$$
M\cdot {g \over h}-{f \over h},\qquad {f \over h}-m \cdot{g \over h}
\leqno(3.5)
$$
are synchronous. Then we have the inequalities
$$
\aligned
0& \le \left| \matrix \int_\Omega \phi f^2 & \int_\Omega \phi f h\\
                      \int_\Omega \phi f h & \int_\Omega \phi h^2\\
              \endmatrix\right|
\cdot
       \left|\matrix \int_\Omega \phi g^2 & \int_\Omega \phi g h\\
                      \int_\Omega \phi g h & \int_\Omega \phi h^2\\
              \endmatrix\right|\\
&\quad - \left|\matrix \int_\Omega \phi f g & \int_\Omega \phi f h\\
                \int_\Omega \phi g h & \int_\Omega \phi h^2\\
        \endmatrix\right|^2\\
&\le {1 \over 4}(M-m)^2
    \left|\matrix \int_\Omega \phi g^2 &  \int_\Omega \phi g h\\
                  \int_\Omega \phi g h & \int_\Omega \phi h^2\\
              \endmatrix\right|^2.
\endaligned
\leqno(3.6)
$$
The constant  ${1 \over 4}$ is best possible in $(3.6)$.}
\vskip 2mm

The proof is obvious by Theorem 2.1 and we omit the details.
\vskip 0.2cm 

{\bf Proposition 3.2.} {\it With the assumptions of Proposition 3.1, we have the inequality
$$\aligned
0& \le \left| \matrix \int_\Omega \phi f^2 & \int_\Omega \phi f h\\
                      \int_\Omega \phi f h & \int_\Omega \phi h^2\\
              \endmatrix\right|^{1 \over 2}
       \left| \matrix \int_\Omega \phi g^2 & \int_\Omega \phi g h\\
                      \int_\Omega \phi g h & \int_\Omega \phi h^2\\
              \endmatrix\right|^{1 \over 2}\\
&\le {1 \over 2}\cdot{(M- m) \over \sqrt{mM}}
    \left| \matrix \int_\Omega \phi f g & \int_\Omega \phi f h\\
                      \int_\Omega \phi g h & \int_\Omega \phi h^2\\
              \endmatrix\right|.
\endaligned
\leqno(3.7)
$$
The constant  ${1 \over 2}$ is best possible in $(3.7)$.}
\vskip 2mm

The following counterpart of the (CBS)-inequality for determinants also holds. 
\vskip 2mm

{\bf Proposition 3.3.} {\it With the assumptions of Proposition 3.1, we have the inequalities
$$
\aligned
0& \le \left| \matrix \int_\Omega \phi f^2 & \int_\Omega \phi f h\\
                      \int_\Omega \phi f h & \int_\Omega \phi h^2\\
              \endmatrix\right|^{1 \over 2}
       \left| \matrix \int_\Omega \phi g^2 & \int_\Omega \phi g h\\
                      \int_\Omega \phi g h & \int_\Omega \phi h^2\\
              \endmatrix\right|^{1 \over 2}\\
&\quad - \left| \matrix \int_\Omega \phi f g & \int_\Omega \phi f h\\
                      \int_\Omega \phi g h & \int_\Omega \phi h^2\\
              \endmatrix\right|\\
&\le {1 \over 2}{({\sqrt M-\sqrt m})^2 \over \sqrt{mM}}
    \left| \matrix \int_\Omega \phi f g & \int_\Omega \phi f h\\
                      \int_\Omega \phi g h & \int_\Omega \phi h^2\\
              \endmatrix\right|
\endaligned
\leqno(3.8)
$$
and
$$
\aligned
0& \le \left| \matrix \int_\Omega \phi f^2 & \int_\Omega \phi f h\\
                      \int_\Omega \phi f h & \int_\Omega \phi h^2\\
              \endmatrix\right|
       \left| \matrix \int_\Omega \phi g^2 & \int_\Omega \phi g h\\
                      \int_\Omega \phi g h & \int_\Omega \phi h^2\\
              \endmatrix\right|\\
&\quad- \left| \matrix \int_\Omega \phi f g & \int_\Omega \phi f h\\
                 \int_\Omega \phi g h & \int_\Omega \phi h^2\\
              \endmatrix\right|^2\\
&\le {1 \over 4}{(M-m)^2 \over \sqrt{mM}}
    \left| \matrix \int_\Omega \phi f g & \int_\Omega \phi f h\\
                   \int_\Omega \phi g h & \int_\Omega \phi h^2\\
              \endmatrix\right|^2.
\endaligned
\leqno(3.9)
$$
The constants ${1 \over 2}$ in $(3.8)$ and ${1 \over 4}$ in $(3.9)$ are best possible.}
\vskip 0.2cm 

Finally, we may state the following reverse of the triangle inequality for determinants:
\vskip 2mm

{\bf Proposition 3.4.} {\it With the assumptions of Proposition 3.1, we have the inequalities
$$
\aligned
0& \le \left|\matrix \int_\Omega \phi f^2 & \int_\Omega \phi f h\\
                      \int_\Omega \phi f h & \int_\Omega \phi h^2\\
              \endmatrix\right|^{1 \over 2}
 + \left|\matrix \int_\Omega \phi g^2 & \int_\Omega \phi g h\\
                      \int_\Omega \phi g h & \int_\Omega \phi h^2\\
              \endmatrix\right|^{1 \over 2}\\
& - \left|\matrix \int_\Omega \phi (f+g)^2 & \int_\Omega \phi (f+g) h\\
                      \int_\Omega \phi (f+g) h & \int_\Omega \phi h^2\\
              \endmatrix\right|^{1 \over 2}\\
&\le {\sqrt M-\sqrt m) \over (mM)^{1 \over 4}}
    \left|\matrix \int_\Omega \phi f g & \int_\Omega \phi f h\\
                      \int_\Omega \phi g h & \int_\Omega \phi h^2\\
              \endmatrix\right|^{1 \over 2}.
\endaligned
\leqno(3.10)
$$}
\vskip 0.2cm

{\bf Acknowledgement.} S. S. Dragomir and Y. J. Cho greatly acknowledge the financial support from Brain Pool
Program (2002) of the Korean Federation of Science and Technology Societies. The research was performed under the "Memorandum of 
Understanding" between Victoria University and Gyeongsang National University.
\vskip 0.5cm

\centerline{\bf References} 
\vskip 0.2cm
\item{ [1]} C. Budimir, Y. J. Cho, M. Mati\'c and J. Pe\v cari\' c, Cebysev's inequality in n-inner product space,
``{\it Inequality Theory and Applications}, Vol. 2" (Editors: Y. J. Cho, J. K. Kim and S. S. Dragomir), Nova Science Publishers, Inc.,  New York, 2001, 
pp. 87--94.
\item{\ [2]} Y. J. Cho, S. S. Dragomir, A. White and S. S. Kim, Some inequalities in 2-inner product spaces,
{\it Demonstratio Math.} {\bf 32}(3)(1999),  485--493.
\item{\ [3]} Y. J. Cho, Paul C. S Lin, S. S. Kim and A. Misiak, {\it Theory  of 2-Inner Product Spaces},
 Nova Science Publishers, Inc.,  New York, 2001.
\item{\ [4]} Y. J. Cho, M. Mati\' c and J. Pe\v cari\' c, On Gram's determinant in 2-inner product spaces,
 {\it J. Korean Math. Soc.}, {\bf 38}(6)(2001),  1125--1156.
\item{\ [5]} Y. J. Cho, M. Mati\' c and J. Pe\v cari\' c, Inequalities of Hlawka's type in G-inner product spaces,
 ``{\it Inequality Theory and Applications}, Vol. 1" (Editors:  Y. J. Cho, J. K. Kim and S. S. Dragomir),
 Nova Science Publishers, Inc., New York, 2001, pp. 95--102.
\item{\ [6]} Y. J. Cho, M. Mati\'c and J. Pe\v cari\' c, On Gram's determinant in n-inner product spaces,
preprint.
\item{\ [7]} S. S. Dragomir, Y. J. Cho and S. S. Kim, Superadditivity and monotonicity of 2-norms generated
by inner products and related results, {\it Soochow. J. Math.} {\bf 29}(1)(1998), 13--32.
\item{\ [8]} R. W. Freese and Y. J. Cho, {\it Geometry of Linear 2-Normed Spaces}, Nova Science Publishers, Inc., New York, 2001.
\item{\ [9]} H. Gunawan, On n-inner product, n-norms and the Cauchy-Schwarz inequality,
 {\it Sci. Math. Japon.} {\bf 55}(1)(2002),  53--60.
\item{\ [10]} C. S. Lin, On inequalities in inner product and 2-inner product spaces, {\it Internat. J. Pure and Appl. Math.} {\bf 3}(3)(2002), 287--298.
\item{[11]} S. S. Kim, S. S. Dragomir, A. White and Y. J. Cho, On the Gr\" uss type inequality in
 2-inner product spaces and applications, {\it PanAmer. Math. J.} {\bf 71}(3)(2001),  89--97.
\item{[12]} S. S. Kim and S. S. Dragomir, Inequalities involuding Gram's determinant in 2-inner product spaces,
``{\it Inequality Theory and Applications}, Vol. 1" (Editors: Y. J. Cho, J. K. Kim and S. S. Dragomir), Nova Science Publishers, Inc., New York, 2001,
pp.183--192. 
\vskip 5mm

Sever S. Dragomir\par
School of Computer Science and Mathematics\par
Victoria University\par
P. O. Box 14428, Melbourne City MC\par
Victoria 8001, Australia\par
{\it E-mail: sever.drgomir\@vu.edu.au}
\vskip 2mm

Yeol Je Cho\par
Department of Mathematics Education\par
The Research Institute of Natural Sciences\par
College of Education\par
Gyeongsang National University\par
Chinju 660-701, Korea\par
{\it E-mail: yjcho\@nongae.gsnu.ac.kr}
\vskip 2mm

Seong Sik Kim\par
Department of Mathematics\par
Dongeui University\par
Pusan  614-714, Korea\par
{\it sskim\@dongeui.ac.kr}

\end